%% file: paper.tex
\ifdefined\MPC
  \RequirePackage{fix-cm}
  \documentclass[smallextended]{svjour3}
  \smartqed
\else
  \documentclass[11pt,a4paper]{article}
  \usepackage[margin=25mm]{geometry}
  \usepackage{authblk}  % superscript-numbered affiliations, sn-jnl style

  \newcommand{\keywords}[1]{}
  \newcommand{\subclass}[1]{}
  \newcommand{\ackname}{Acknowledgements}
  \newenvironment{acknowledgements}%
    {\par\addvspace{6pt plus 1ex minus .2ex}\noindent
     \textbf{\ackname}\hskip 1em\ignorespaces}%
    {\par\addvspace{6pt}}
\fi

\usepackage{graphicx}
\usepackage{tikz}
\usetikzlibrary{arrows.meta,positioning,fit,backgrounds,calc}
\usepackage{pgfplots}
\pgfplotsset{compat=1.18}
\usepackage{amsmath,amssymb}
\usepackage{booktabs}
\usepackage{algorithmic}
\usepackage{algorithm}
\usepackage{hyperref}
\usepackage{xcolor}
\usepackage{xspace}
\usepackage{todonotes}

\newcommand{\scip}{\textsc{SCIP}\xspace}

\newcommand{\fiberscip}{\textsc{FiberSCIP}\xspace}
\newcommand{\rexi}{\textsc{ReXi}\xspace}
\newcommand{\rexils}{\textsc{ReXiLS}\xspace}
\newcommand{\rexilns}{\textsc{ReXiLNS}\xspace}
\begin{document}

\title{Race, Exchange, Improve: Finding high-quality MIP solutions quickly}

\ifdefined\MPC
\journalname{Mathematical Programming Computation}

\author{Gioni Mexi \and Daniel Rehfeldt}

\institute{
	G.~Mexi \and D.~Rehfeldt \at
	Zuse Institute Berlin, Takustra\ss e 7, 14195 Berlin, Germany \\
	\email{\{mexi, rehfeldt\}@zib.de}
	\and
	D.~Rehfeldt \at
	IVU Traffic Technologies AG, Bundesallee 88, 12161 Berlin, Germany
}

\date{Received: date / Accepted: date}
\else
\author[1]{Gioni Mexi}
\author[1,2]{Daniel Rehfeldt}
\affil[1]{Zuse Institute Berlin, Takustra\ss e 7, 14195 Berlin, Germany}
\affil[2]{IVU Traffic Technologies AG, Bundesallee 88, 12161 Berlin, Germany}
\affil[ ]{\texttt{\{mexi, rehfeldt\}@zib.de}}
\date{}
\fi

\maketitle

\begin{abstract}
	
	Mixed-integer programming (MIP) is a cornerstone in applied optimization, both in industry and academia.
	Recently, there has been increased attention to finding strong primal solutions quickly. This is reflected, for example, in the development of the NVIDIA cuOpt solver and, most recently, in the new MIPFEAS benchmark, which has a tight time limit of 600 seconds and evaluates solvers based on how quickly they find high-quality primal solutions.
	
	This article introduces a MIP portfolio parallelization scheme, focusing on efficiently exchanging information between its workers. We present two implementations of this scheme: one built directly into the open-source MIP solver \scip, and an external one, which we call \rexi.
   \rexi is currently the fastest non-commercial solver in the MIPFEAS benchmark, followed by the \scip-integrated implementation. Moreover, we present new versions of both implementations that
   considerably outperform their predecessors on the MIPFEAS benchmark.  

\keywords{Mixed-integer programming \and Parallel solving \and Algorithm racing}
\subclass{90C11 \and 90C06 \and 68W10}
\end{abstract}

%----------------------------------------------------------------------
\section{Introduction}
\label{sec:intro}
%----------------------------------------------------------------------

Mixed-integer programming solvers have long been a standard tool for solving optimization problems in virtually all industry sectors as well as in academia. Several powerful commercial (e.g., CPLEX~\cite{cplex125}, Gurobi~\cite{gurobi}, Xpress~\cite{BelottiBertholdGallyGottwaldPolik2025}, and COPT~\cite{ge2022cardinal}) and academic (e.g., \scip~\cite{hojny2025scip}, HiGHS~\cite{huangfu2018parallelizing}, and CBC~\cite{forrest2005cbc}) solvers are available. Recently, there has been increased attention to finding strong primal solutions quickly. Examples are the development of the NVIDIA cuOpt solver~\cite{nvidia2025cuopt} and, most recently, the release of the MIPFEAS benchmark~\cite{mipfeas}, which has a tight time limit of 600 seconds and evaluates solvers based on the primal integral~\cite{BERTHOLD2013611}, which favors methods that quickly find high-quality solutions. At the same time, improved hardware, such as GPUs and modern multicore CPUs, provides additional opportunities for accelerating MIP solving. 
Running multiple solver configurations in parallel can exploit their complementary strengths, an idea already used in parallel solver frameworks such as \fiberscip~\cite{shinano2018fiberscip} and CP-SAT~\cite{perron_et_al:LIPIcs.CP.2023.3}. An early demonstration for MIP is due to Carvajal et al.~\cite{carvajal2014using}, who race diversified workers that exchange incumbent solutions and bounds.

This article introduces a shared-memory parallel framework, based on \scip, aimed at finding high-quality solutions within short running times. It is implemented both directly in \scip, extending its concurrent solving framework, and in a new, more specialized solver called \rexi.
Both implementations build a racing
portfolio of diversified \scip workers, each using one thread. No coordinator thread is used,
but information is shared between workers through a solution pool and a bound pool, both designed for
minimal blocking. Importantly, some of the workers are run without presolving, which allows finding solutions more quickly, while some of the lost problem strengthening is regained via the bound pool (which communicates some of the presolving results of other workers).    
With the focus being primal feasibility, \rexi devotes two of its workers
entirely to primal heuristics, a highly optimized local search and a large
neighborhood search,
rather than to branch-and-bound search. 
We report a computational study that separates the
contribution of racing, of each exchange mechanism, and of the two primal
workers.
While none of the components is completely new, the main contribution of this article lies
in incorporating them into a lean design, based on \scip, with a carefully
optimized implementation. Using ten threads, this improves the primal integral by
more than a factor of five over single-threaded \scip.

\rexi is currently the fastest non-commercial solver in the MIPFEAS benchmark~\cite{mipfeas}, followed by concurrent \scip. Moreover, the latest version of \rexi is considerably faster and is on par with the average performance of three leading commercial solvers on this benchmark. It should be noted, however, that these commercial solvers are run with their default settings, whereas \rexi is specifically designed to find high-quality primal solutions quickly.

The remainder of this article is organized as follows.
Section~\ref{sec:racing} describes the framework along its three pillars: the
racing portfolio, the exchange of solutions and bounds, and the dedicated primal
workers. Section~\ref{sec:impl} presents the two implementations, and where they
differ. Finally, Section~\ref{sec:results} reports the computational study.

%----------------------------------------------------------------------
\section{Race, Exchange, Improve}
\label{sec:racing}
%----------------------------------------------------------------------

The architecture is based on three pillars, encoded in the name of the new solver \rexi\footnote{Which, by sheer luck, also happens to fit the last names of the authors.}:

\begin{itemize}
  \item \emph{Race}: MIP solvers are well known for their performance variability~\cite{Lodi_2013}. A diversified portfolio can transform this behavior into an asset~\cite{huberman1997,gomes2001,fischetti2014erraticism}.
  \item \emph{Exchange}: We aim to share as much relevant information as possible among workers, while keeping overhead minimal. This includes the exchange of primal solutions, primal bounds, and variable bounds.
  \item  \emph{Improve}: LP-free heuristics can sometimes outperform MIP workers. Thus we include such heuristics in our solver. Instead of calling them within a MIP worker, we give them a full thread for continuous work, but still tightly integrate them with the MIP workers via mutual solution and primal bound exchange.
\end{itemize}

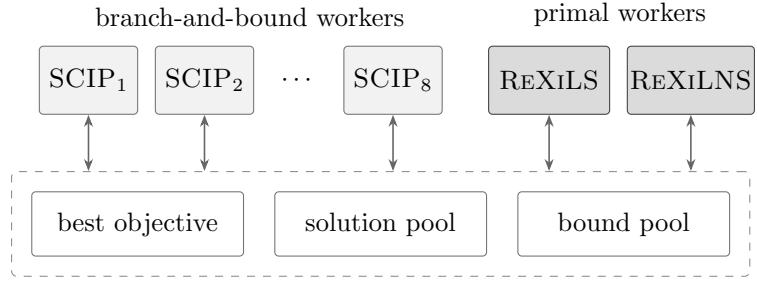
\begin{figure}[t]
  \centering
  \begin{tikzpicture}[
    font=\small,
    bb/.style={draw=black!55, fill=black!5, rounded corners=1.5pt,
               minimum width=13mm, minimum height=9mm, align=center,
               inner sep=1pt, font=\small},
    side/.style={bb, draw=black!75, fill=black!14, minimum width=16mm},
    store/.style={draw=black!55, fill=white, rounded corners=1.5pt,
                  minimum height=8.5mm, minimum width=28mm, inner xsep=2.5mm,
                  align=center, font=\small},
    lnk/.style={{Stealth[length=1.6mm]}-{Stealth[length=1.6mm]},
                draw=black!60, semithick},
  ]
    % ---- workers ---------------------------------------------------------
    \node[bb]                     (w1) at (0,0)      {\textsc{SCIP}$_1$};
    \node[bb, right=2.2mm of w1]  (w2)               {\textsc{SCIP}$_2$};
    \node[right=2.2mm of w2, font=\small]  (dots) {$\cdots$};
    \node[bb, right=2.2mm of dots](w8)               {\textsc{SCIP}$_8$};
    \node[side, right=6mm of w8]  (ls)               {\rexils};
    \node[side, right=2.2mm of ls](lns)              {\rexilns};

    \node[font=\small, above=1.2mm of w2.north, xshift=6mm]
      {branch-and-bound workers};
    \node[font=\small, above=1.2mm of $(ls.north)!0.5!(lns.north)$]
      {primal workers};

    % ---- shared layer ----------------------------------------------------
    \node[store] (obj) at (6.5mm,-19mm) {best objective};
    \node[store, right=4mm of obj] (sol) {solution pool};
    \node[store, right=4mm of sol] (bnd) {bound pool};

    \begin{scope}[on background layer]
      \node[draw=black!45, dashed, rounded corners=2pt, inner sep=2.6mm,
            fit=(obj)(sol)(bnd)] (layer) {};
    \end{scope}

    % ---- links -----------------------------------------------------------
    \foreach \w in {w1,w2,w8,ls,lns} {\draw[lnk] (\w.south) -- (\w.south |- layer.north);}

  \end{tikzpicture}
  \caption{Architecture of \rexi. Every worker races the same model and
    communicates only through the shared structures (Section~\ref{sec:ex}). Eight workers run branch-and-bound
    with diversified settings (Section~\ref{sec:race}); \rexils and \rexilns own
    a thread each (Section~\ref{sec:impr}). Concurrent
    \scip uses the same layer with ten branch-and-bound workers.}
  \label{fig:architecture}
\end{figure}

Much effort went into spending as little time as possible on exchanging
information. Notably, we do not use a coordinator or supervisor thread, but exchange information via shared data structures. The communication is not lock-free, but in practice has minimal contention and very short waiting times.
Figure~\ref{fig:architecture} gives an overview.

%----------------------------------
\subsection{Race: The Portfolio}
\label{sec:race}
%----------------------------------

\rexi races $N$ independent \scip instances on the same model, each with
its own parameter setting and its own random seed. We use the following ten
\scip settings: default \scip; no separation combined with fast presolve;
two emphasis \emph{feasibility} settings that differ only through seed
diversification; no presolving; SAT-like depth-first search without any LP
solves; periodic feasibility jump; periodic local search; a local search running
perpetually; and depth-first search with aggressive restarting and primal
heuristics.

These ten settings are exactly the portfolio raced by \scip's concurrent mode
and \rexi. Note that \rexi runs \rexils perpetually instead of \scip's local
search, and replaces another setting by a large neighborhood search (LNS) worker, as described in
Section~\ref{sec:impr}.

The workers do different amounts of presolving, some none at all, and use
different search strategies, primarily to find good solutions quickly, but also
to diversify the search.

%----------------------------------
\subsection{Exchange: Solutions and Bounds}
\label{sec:ex}
%----------------------------------

We mostly exchange two classes of solving information among the workers.
First, primal solutions and primal bounds, and, second, variable bounds.

For solution and primal bound exchange, we use a two-lane approach:
The fast lane exchanges only the current best primal bound. The exchange has negligible overhead but still delivers very important information (for example, to cut off branch-and-bound trees or for reduced-cost fixing).
The slow lane exchanges actual primal solutions via a pool described in the following. The solutions are exchanged less often to minimize overhead. More accurately: they are pulled less often from the pool, but published right away. The latter is especially important for the two side-heuristics described in the next subsection.
For the solution pool, we use the following:
\begin{itemize}
	\item Solutions in the pool live in the original (unpresolved) space. This is important because workers do different presolving.
	\item Each worker publishes a new solution only if it is better than the current best one in the pool.
	\item Each publisher is responsible for checking the feasibility of each solution in the original space before publishing.
	\item We avoid waiting times when publishing solutions by using a pointer-based pool. Each worker creates the best solution in the original space, then checks for feasibility, and only then locks the pool for simply adding a pointer to this new solution.
     \item As to pulling solutions (which is done less often): we transform the current best solution to the presolved space, which might lead to an invalid solution. In this case, we continue with the next best solution from the pool and so on.
\end{itemize}

In contrast to the solution pool, for bound exchange we use fixed-size data structures to store globally valid variable bounds, which are strengthened through the solution process. Such sharing of bound tightenings between concurrently solving workers goes back to distributed domain propagation~\cite{gottwald2017distributed}.
\begin{itemize}
\item Workers share bounds from domain propagation, but also from presolving. The latter benefits especially the workers that were run without presolving.
\item Improved bounds are published immediately. New bounds are only pulled at specific points, since the additional domain propagations that are triggered in \scip after bound changes can be expensive.
\item Strong dual reductions can only be performed by one worker. Otherwise, we could potentially cut off all optimal solutions.
\end{itemize}
The main weak points of the bound pool used for \scip are that strong dual reductions are confined to one worker, plus the cost of \scip's re-propagation of new bounds.

More details on the implementation of both the bound and solution pool can be found in Section~\ref{sec:impl}.

%----------------------------------
\subsection{Improve: Primal Workers}
\label{sec:impr}
%----------------------------------

There are three types of workers:

\begin{itemize}
	\item Standard branch-and-bound \scip workers. These also work on the dual side, which in turn allows for finding better primal solutions.
	\item One pure primal LP-free heuristic that runs perpetually.
	\item One improvement heuristic that perpetually uses solutions added by other workers to the common pool and recombines them, running LNS.
\end{itemize}

The local search worker runs an LP-free local search~\cite{LIN2025104405} on a minimally presolved problem. This heuristic maintains a single complete variable assignment and improves it by single-variable moves scored on constraint violation and objective value.
The heuristic is restarted from the best incumbent in the pool whenever
its own search stalls. This heuristic is especially important for instances that spend a significant amount of time in presolving or the initial root LP, because \scip runs most other primal heuristics only afterwards.

The LNS~\cite{shaw1998,danna2005rins} worker needs at least one reference point to define the search neighborhood. It takes the best solution in the shared pool; 
while the pool is still empty it takes the latest retained point of smallest constraint violation from the local search. 

From that reference it builds a neighborhood in one of two ways.
Mutation fixes each integer variable independently with probability $\alpha$ to its rounded reference value.
Crossover fixes the agreement set of two pooled solutions from different workers. Whenever an incumbent exists it is also imposed as an objective limit, so a
neighborhood that contains nothing better is proven empty and abandoned quickly.

The crossover mode is closely related to the evolutionary algorithm for polishing
MIP solutions of Rothberg~\cite{rothberg2007evolutionary}, known as solution
polishing in CPLEX. Dedicating whole workers to large neighborhood search over a
shared solution pool is also done by CP-SAT~\cite{perron_et_al:LIPIcs.CP.2023.3}.

\section{Two implementations: Concurrent SCIP and ReXi}
\label{sec:impl}
%----------------------------------------------------------------------

This section describes the two implementations of the race-exchange-improve framework introduced in the previous section. Both are based on \scip, but while the first one is directly implemented in \scip, the second, \rexi, implements several key components outside of \scip for increased efficiency. In this way, one can use more specialized data structures and implementations.

\subsection{Concurrent SCIP}

Our implementation builds on the existing concurrent \scip framework, introduced in \scip~4.0~\cite{MaherFischerGallyetal.2017}. 
That framework already races a portfolio of \scip workers, with the difference that information is exchanged at fixed synchronization points.

A solver that reaches the synchronization point writes its solutions and bounds into a shared store and then waits until every other solver has reached that point before reading. Therefore, the exchange is only as fast as the slowest worker, and a solver that has just found a good solution cannot pass it on until the others have caught up.

The solution and bound pools of Section~\ref{sec:ex} avoid this.
We use a fast size-hint atomic to check for each consumer whether anything has been added to the pool since the last pull. 
The bound pool is relatively simple: static arrays, with the whole pool locked for each update or pull.
Concurrent \scip also uses a local search heuristic~\cite{LIN2025104405} that runs
alongside the solve, which will be part of the \scip~11 release.

\subsection{ReXi}

\rexi (including its two primal workers) is written in C++, making use of C++ parallelization utilities, such as threads, atomics and mutexes. 

The solution pool is mostly implemented as described in Section~\ref{sec:ex}. As in concurrent \scip, we use a size-hint atomic to check before each pull whether anything new has been added to the pool. The best primal bound is shared via a single atomic variable.

The bound pool implementation includes some additional details compared to those given in Section~\ref{sec:ex}. We experimented with bound pools from the literature, including the ``dirty sets'' implementation of OR-Tools/CP-SAT. However, in our application, we have the somewhat special case that, on the one hand, we have very aggressive updates of bound changes at the beginning of the solve, because we update the bounds also during presolving. On the other hand, different workers pull these bounds aggressively. On some instances this leads to noticeable blocking times (even though it never costs more than a few percent of runtime).
We use the following design:

\begin{itemize}
\item Three contiguous, fixed-size arrays over the original variables: lb, ub, var\_version.  Bounds are
  monotone and overwritten in place, so memory is constant regardless
  of the number of tightenings, and a worker restarting its search can
  recover the complete current state with a single pull (important for workers without presolving).
\item The variable range is split into 32 contiguous shards, each with its own mutex and its own atomic published version. Updating or reading bounds locks exactly one shard, leading to reduced contention. 
\item Per-reader state (seen versions, counters); the reader-owned fields are touched only by the owning thread.
\item Change detection is a three-level mechanism, to minimize blocking. As in the solution pool,
  one atomic global counter gives a lock-free fast path (to signify if anything has
  changed at all); the per-shard counters let a reader skip untouched
  shards without locking; and the per-variable versions give exactly
  the changed bounds within a locked shard.
\end{itemize}

For each instance of the MIPFEAS benchmark, the waiting time of the above bound pool is far below 1 percent of the overall runtime, so completely negligible. Still, as already mentioned, there is a significant performance penalty from the rather slow repropagations triggered in \scip. Similarly, using strong dual reductions only on one worker (as required when the bound pool is active) costs performance.

For the local-search worker \rexils we mostly follow~\cite{LIN2025104405}, but with an optimized implementation that is significantly faster than the one from~\cite{LIN2025104405}. For some key instances, for which local search finds the first good solution, we observe a speed-up of more than a factor of three compared to the original implementation.

Besides \rexilns, described already in Section~\ref{sec:impr}, another difference of \rexi compared to concurrent \scip is the solution feasibility checking which happens before each publication in the solution pool. In \rexi this is performed with a more cache-efficient CSR-based check, which still uses the same tolerances as native \scip, but is faster.

%----------------------------------------------------------------------
\section{Computational Study}
\label{sec:results}
%----------------------------------------------------------------------
In this section, we measure the performance gains from racing, information sharing, the implementation differences and the two specialized workers of \rexi. 
Experiments are run on AMD EPYC 9B45 CPUs with up to $10$ threads, a
$96$\,GB memory limit and a $600$\,s time limit, with cluster nodes used exclusively.
All runs use pre-release versions of \scip~11 and its LP solver SoPlex~9.
Our testset consists of the $233$ instances of the MIPLIB 2017
benchmark~\cite{gleixner2021miplib}, excluding instances known to be infeasible.

For primal feasibility, the metric we use is the primal integral, in the form used
by the MIPFEAS benchmark~\cite{mipfeas}, which adapts the original definition of
Berthold~\cite{BERTHOLD2013611}. Let $\bar{x}(t)$ be the objective value of
the incumbent at time $t$ and $x^\star$ the reference value, and let
\[
  p(t) =
  \begin{cases}
    2, & \text{no feasible solution is known at time } t, \\
    1, & \bar{x}(t) \text{ and } x^\star \text{ have opposite signs}, \\
    0, & |\bar{x}(t)| < 10^{-6} \text{ and } |x^\star| < 10^{-6}, \\[2pt]
    \dfrac{|\bar{x}(t) - x^\star|}{\max\{|\bar{x}(t)|,\,|x^\star|\}}, & \text{otherwise}
  \end{cases}
\]
be the penalty incurred at time $t$. Since $p$ changes only when a new incumbent
is found, it is piecewise constant, and its average over the time limit $T$ is a
finite sum: with $0 = t_0 < t_1 < \dots < t_k \le T$ the times at which the
incumbent improves and $t_{k+1} = T$,
\[
  P(T) \;=\; \frac{1}{T} \sum_{i=0}^{k} p(t_i) \, (t_{i+1} - t_i) \;\in\; [0,2] .
\]
A value of $2$ means that no feasible solution was found within the time limit;
smaller values mean that good solutions were found earlier.

In the tables below, ``opt.\ found'' counts instances whose reference objective
was reached and ``opt.\ proven'' those on which the run closed its own
primal--dual gap. Time and primal integral (PI) are shifted geometric means with
shifts $1$\,s and $0.001$.

Table~\ref{tab:scaling} gives the overall picture. It compares default \scip,
\scip's concurrent mode, and \rexi. The row marked ``concurrent settings'' races the
same ten settings as \scip's concurrent mode (Section~\ref{sec:race}); the last
row is \rexi in its default configuration, which replaces one of those settings
by the LNS worker.

\begin{table}[ht]
  \centering
  \caption{Performance comparison of default \scip, \scip's concurrent mode, and
    \rexi with and without its LNS worker.}
  \label{tab:scaling}
  \small
  \setlength{\tabcolsep}{4pt}
  \begin{tabular}{@{}lrrrrrr@{}}
    \toprule
    & \multicolumn{1}{c}{threads}
    & \multicolumn{1}{c}{feas.}
    & \multicolumn{1}{c}{opt.\ found}
    & \multicolumn{1}{c}{opt.\ proven}
    & \multicolumn{1}{c}{time [s]}
    & \multicolumn{1}{c}{PI} \\
    \midrule
    \scip                       & $1$  & $211$ & $119$ & $95$ & $201.2$ & $0.0568$ \\
    \scip concurrent            & $10$ & $224$ & $154$ & $119$ & $149.6$ & $0.0178$ \\
    \rexi (concurrent settings) & $10$ & $\mathbf{225}$ & $158$ & $123$ & $142.8$ & $0.0134$ \\
    \midrule
    \rexi                       & $10$ & $\mathbf{225}$ & $\mathbf{160}$ & $\mathbf{124}$
                                  & $\mathbf{135.3}$ & $\mathbf{0.0107}$ \\
    \bottomrule
  \end{tabular}
\end{table}

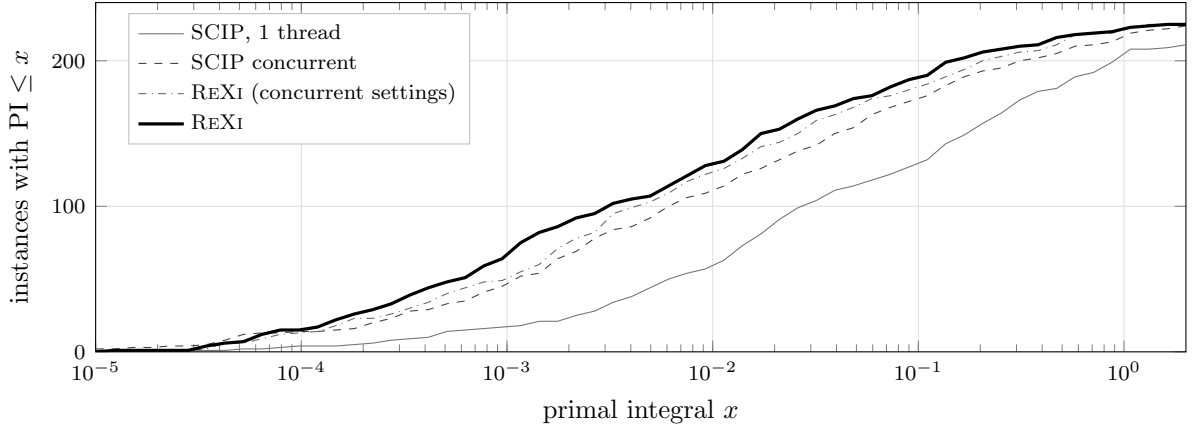
\begin{figure}[ht]
  \centering
  \input{data/pi_profile}
  \caption{Cumulative distribution of the primal integral over the $233$
    instances, for the four configurations of Table~\ref{tab:scaling}. For a
    threshold $x$ the curve gives the number of instances reaching a primal
    integral of at most $x$.}
  \label{fig:pi-profile}
\end{figure}

Going from one thread to ten cuts the shifted geometric mean of the primal
integral by $69\%$, from $0.0568$ to $0.0178$. It also adds $13$ instances on
which a feasible solution is found at all, and $35$ on which the reference
optimum is reached. Comparing the two \scip runs against each other instance by
instance, the ten-thread configuration is at least $10\%$ better on the primal
integral on $202$ of the $233$ instances, against $7$ for the single-thread run.
This is not the effect of a small subset of instances. As
Figure~\ref{fig:pi-profile} shows, the two distributions are separated over the
whole range.

The last two rows show \rexi. The third row races the same settings as
concurrent \scip, so the two differ only in the engine that runs the race and
in the \rexils implementation. The shifted geometric mean
of the primal integral is $0.0134$ against $0.0178$.
The solution pool then enables further ideas, such as the LNS worker of
Section~\ref{sec:impr}, which in the last row replaces one of the portfolio
settings. This brings the primal integral to $0.0107$, another $20\%$ below the
row above it and $40\%$ below concurrent \scip, cuts solving time by $10\%$
against concurrent \scip and reaches the reference optimum on six more
instances. We see the same effect on the dual side, where \rexi proves the most
optima of the four, $124$ against concurrent \scip's $119$.

Table~\ref{tab:exchange} switches off one of the two exchange mechanisms at a time, leaving the ten-thread portfolio otherwise unchanged.
Note that without solution exchange the LNS worker has nothing to build
neighborhoods from, so it is not used and its seat falls back to the default
portfolio setting, i.e.\ to the \rexi configuration of the third row of
Table~\ref{tab:scaling}.
\begin{table}[ht]
  \centering
  \caption{Effect of disabling one of the two exchange mechanisms in \rexi, ten threads.}
  \label{tab:exchange}
  \small
  \setlength{\tabcolsep}{4pt}
  \begin{tabular}{@{}lrrrrr@{}}
    \toprule
    & \multicolumn{1}{c}{feas.}
    & \multicolumn{1}{c}{opt.\ found}
    & \multicolumn{1}{c}{opt.\ proven}
    & \multicolumn{1}{c}{time [s]}
    & \multicolumn{1}{c}{PI} \\
    \midrule
    \rexi                     & $\mathbf{225}$ & $\mathbf{160}$ & $\mathbf{124}$
                                & $\mathbf{135.3}$ & $\mathbf{0.0107}$ \\
    \quad no solution exchange  & $224$ & $150$ & $115$ & $156.0$ & $0.0146$ \\
    \quad no bound exchange     & $224$ & $157$ & $\mathbf{124}$ & $137.6$ & $0.0118$ \\
    \bottomrule
  \end{tabular}
\end{table}

Of the two mechanisms, solution exchange turns out to be the more important
one. Without it the shifted
geometric mean of the primal integral rises by $37\%$, from $0.0107$ to
$0.0146$, solving time by $15\%$, and the number of instances solved to proven
optimality falls from $124$ to $115$. The workers also find around $22\%$ more
solutions in total, since each of them has to rediscover solutions that another
worker might have already discovered.

Switching off bound exchange costs $10\%$ on the shifted geometric mean, $0.0107$
to $0.0118$, one instance on which a feasible solution is found and three on which
the reference optimum is reached, while the number of instances closed to proven
optimality is unchanged at $124$ and solving time rises by less than $2\%$. So
the shared bounds help the workers find better solutions sooner, but they do not
increase the number of instances on which we can prove optimality.

Both concurrent \scip and \rexi have been part of the MIPFEAS
benchmark~\cite{mipfeas}. Table~\ref{tab:version} compares the configuration
above against the \rexi build used there. The improvements made since, including
\rexils, bound sharing and the LNS worker, have lowered the shifted geometric mean
of the primal integral by $35\%$, from
$0.0164$ to $0.0107$; the reference optimum is reached on $12$ more instances and
proven on $8$ more, and the shifted geometric mean of the solving time falls by
$11\%$. Per instance the current build is at least $10\%$ better on $165$ of the
$233$ instances, against $26$ the other way. On the benchmark itself the submitted
build scores $0.0165$ against $0.0105$ for the virtual mean commercial
solver~\cite{mipfeas}, so the $0.0107$ of the current build (using a cross-machine comparison) is on par with the average 
of three leading commercial solvers.
Default \scip and its concurrent mode have likewise improved since their benchmark
submissions, mainly through the introduction of the local search
heuristic and improvements in the LP solver SoPlex (the latter were motivated by the development of \rexi).

\begin{table}[ht]
  \centering
  \caption{\rexi against the build submitted to the MIPFEAS benchmark, ten threads.}
  \label{tab:version}
  \small
  \setlength{\tabcolsep}{4pt}
  \begin{tabular}{@{}lrrrrr@{}}
    \toprule
    & \multicolumn{1}{c}{feas.}
    & \multicolumn{1}{c}{opt.\ found}
    & \multicolumn{1}{c}{opt.\ proven}
    & \multicolumn{1}{c}{time [s]}
    & \multicolumn{1}{c}{PI} \\
    \midrule
    \rexi                   & $225$ & $\mathbf{160}$ & $\mathbf{124}$
                              & $\mathbf{135.3}$ & $\mathbf{0.0107}$ \\
    \quad MIPFEAS submission  & $225$ & $148$ & $116$ & $151.7$ & $0.0164$ \\
    \bottomrule
  \end{tabular}
\end{table}

%----------------------------------

%----------------------------------------------------------------------
\section{Conclusion}
\label{sec:conclusion}
%----------------------------------------------------------------------

Race, exchange, improve showed significant results. Racing diversified \scip workers, exchanging solutions and variable bounds without a coordinator thread, and spending threads on specialized primal heuristics improves the primal integral on ten threads by more than a factor of five over
sequential \scip. On the MIPFEAS benchmark instances, the latest \rexi is competitive with the average of three leading commercial solvers.

Multiple directions for further research are open. More threads and more specialized workers are among
them, since the portfolio has ten settings today while the machines have many
more cores. Also, workers of a
different kind, such as GPU-based solvers, are left for future research.

%----------------------------------------------------------------------
% Acknowledgments
%----------------------------------------------------------------------
\begin{acknowledgements}
The authors thank the current and former \scip developers, whose work this framework builds upon.
Further thanks go to Yuji Shinano for his work on \fiberscip, which facilitated the use of \scip in this work,
and Michael Bussieck and GAMS for the MIPFEAS benchmarking initiative.
\end{acknowledgements}

\begingroup
\renewcommand{\ackname}{Funding disclosure}
\begin{acknowledgements}
Research reported in this paper was partially supported through the Research Campus Modal
funded by the German Federal Ministry of Education and Research (fund numbers 05M14ZAM,
05M20ZBM) and the Deutsche Forschungsgemeinschaft (DFG) through the DFG Cluster of
Excellence MATH+.
\end{acknowledgements}
\endgroup

%----------------------------------------------------------------------
% References
%----------------------------------------------------------------------
\bibliographystyle{spmpsci}
\bibliography{references}

\end{document}

%% file: data/pi_profile.tex
% generated by data/plot_pi_profile.py -- do not edit by hand
\begin{tikzpicture}
\begin{axis}[
  width=\linewidth, height=6.2cm,
  xmode=log, xlabel={primal integral $x$},
  ylabel={instances with PI $\le x$},
  xmin=1e-05, xmax=2, ymin=0, ymax=240,
  legend pos=north west, legend cell align=left,
  legend style={font=\scriptsize, draw=black!25, fill=white},
  tick label style={font=\scriptsize}, label style={font=\small},
  grid=major, grid style={black!12},
]
\addplot[solid, draw=black!55] coordinates {(1e-05,1) (1.22984e-05,1) (1.5125e-05,1) (1.86013e-05,1) (2.28766e-05,1) (2.81345e-05,1) (3.46009e-05,1) (4.25535e-05,1) (5.2334e-05,2) (6.43623e-05,2) (7.91552e-05,3) (9.73481e-05,4) (0.000119722,4) (0.000147239,4) (0.00018108,5) (0.0002227,6) (0.000273885,8) (0.000336834,9) (0.000414251,10) (0.000509461,14) (0.000626555,15) (0.000770561,16) (0.000947666,17) (0.00116548,18) (0.00143335,21) (0.00176278,21) (0.00216794,25) (0.00266621,28) (0.00327901,34) (0.00403265,38) (0.00495951,44) (0.0060994,50) (0.00750127,54) (0.00922535,57) (0.0113457,63) (0.0139534,73) (0.0171604,81) (0.0211045,91) (0.0259551,99) (0.0319206,104) (0.0392571,111) (0.0482799,114) (0.0593765,118) (0.0730235,122) (0.0898071,127) (0.110448,132) (0.135833,143) (0.167053,149) (0.205448,157) (0.252668,164) (0.310741,173) (0.382161,179) (0.469996,181) (0.578019,189) (0.71087,192) (0.874255,199) (1.07519,208) (1.32231,208) (1.62623,209) (2,211)};
\addlegendentry{\scip, 1 thread}
\addplot[dashed, draw=black!75] coordinates {(1e-05,2) (1.22984e-05,2) (1.5125e-05,3) (1.86013e-05,3) (2.28766e-05,4) (2.81345e-05,4) (3.46009e-05,5) (4.25535e-05,8) (5.2334e-05,12) (6.43623e-05,13) (7.91552e-05,13) (9.73481e-05,14) (0.000119722,14) (0.000147239,15) (0.00018108,16) (0.0002227,20) (0.000273885,23) (0.000336834,28) (0.000414251,29) (0.000509461,33) (0.000626555,35) (0.000770561,41) (0.000947666,45) (0.00116548,52) (0.00143335,54) (0.00176278,64) (0.00216794,69) (0.00266621,78) (0.00327901,84) (0.00403265,86) (0.00495951,92) (0.0060994,100) (0.00750127,106) (0.00922535,109) (0.0113457,114) (0.0139534,122) (0.0171604,126) (0.0211045,132) (0.0259551,138) (0.0319206,142) (0.0392571,150) (0.0482799,154) (0.0593765,163) (0.0730235,168) (0.0898071,172) (0.110448,176) (0.135833,183) (0.167053,189) (0.205448,193) (0.252668,195) (0.310741,200) (0.382161,202) (0.469996,205) (0.578019,210) (0.71087,211) (0.874255,213) (1.07519,219) (1.32231,221) (1.62623,222) (2,224)};
\addlegendentry{\scip concurrent}
\addplot[dash dot, draw=black!60] coordinates {(1e-05,0) (1.22984e-05,1) (1.5125e-05,1) (1.86013e-05,1) (2.28766e-05,1) (2.81345e-05,1) (3.46009e-05,2) (4.25535e-05,6) (5.2334e-05,6) (6.43623e-05,9) (7.91552e-05,12) (9.73481e-05,13) (0.000119722,14) (0.000147239,18) (0.00018108,23) (0.0002227,23) (0.000273885,26) (0.000336834,30) (0.000414251,34) (0.000509461,40) (0.000626555,44) (0.000770561,48) (0.000947666,49) (0.00116548,55) (0.00143335,60) (0.00176278,71) (0.00216794,78) (0.00266621,82) (0.00327901,95) (0.00403265,99) (0.00495951,103) (0.0060994,109) (0.00750127,117) (0.00922535,122) (0.0113457,126) (0.0139534,133) (0.0171604,141) (0.0211045,144) (0.0259551,150) (0.0319206,159) (0.0392571,163) (0.0482799,168) (0.0593765,174) (0.0730235,176) (0.0898071,180) (0.110448,184) (0.135833,189) (0.167053,194) (0.205448,200) (0.252668,203) (0.310741,206) (0.382161,207) (0.469996,211) (0.578019,218) (0.71087,219) (0.874255,219) (1.07519,223) (1.32231,225) (1.62623,225) (2,225)};
\addlegendentry{\rexi (concurrent settings)}
\addplot[solid, very thick, draw=black] coordinates {(1e-05,0) (1.22984e-05,1) (1.5125e-05,1) (1.86013e-05,1) (2.28766e-05,1) (2.81345e-05,1) (3.46009e-05,4) (4.25535e-05,6) (5.2334e-05,7) (6.43623e-05,12) (7.91552e-05,15) (9.73481e-05,15) (0.000119722,17) (0.000147239,22) (0.00018108,26) (0.0002227,29) (0.000273885,33) (0.000336834,39) (0.000414251,44) (0.000509461,48) (0.000626555,51) (0.000770561,59) (0.000947666,64) (0.00116548,75) (0.00143335,82) (0.00176278,86) (0.00216794,92) (0.00266621,95) (0.00327901,102) (0.00403265,105) (0.00495951,107) (0.0060994,114) (0.00750127,121) (0.00922535,128) (0.0113457,131) (0.0139534,139) (0.0171604,150) (0.0211045,153) (0.0259551,160) (0.0319206,166) (0.0392571,169) (0.0482799,174) (0.0593765,176) (0.0730235,182) (0.0898071,187) (0.110448,190) (0.135833,199) (0.167053,202) (0.205448,206) (0.252668,208) (0.310741,210) (0.382161,211) (0.469996,216) (0.578019,218) (0.71087,219) (0.874255,220) (1.07519,223) (1.32231,224) (1.62623,225) (2,225)};
\addlegendentry{\rexi}
\end{axis}
\end{tikzpicture}